\documentclass{article}
\usepackage{amssymb}
\usepackage{graphicx}
\usepackage{amsmath}
\usepackage{epsfig}
\usepackage{amsthm}
\usepackage{url}
\usepackage[all,knot]{xy}
\xyoption{arc}
\newcommand{\cl}[1]{\mathcal{#1}}
\newcommand{\db}[1]{\mathbb{#1}}
\newcommand{\fr}[1]{\mathfrak{#1}}
\newcommand{\fhi}{\varphi}
\newcommand{\fin}{$_{\square}$}
\newcommand{\finf}{_{\square}}
\newcommand{\dis}{\displaystyle}

\newtheorem{theorem}{Theorem}[section]
\newtheorem{lemma}[theorem]{Lemma}

\newtheorem{corollary}[theorem]{Corollary}
\newenvironment{definition}[1][Definition]{\begin{trivlist}
		\item[\hskip \labelsep {\bfseries #1}]}{\end{trivlist}}

\newenvironment{remark}[1][Remark]{\begin{trivlist}
		\item[\hskip \labelsep {\bfseries #1}]}{\end{trivlist}}

\DeclareMathOperator{\ob}{Obj}
\DeclareMathOperator{\ho}{Hom}
\DeclareMathOperator{\mor}{Mor}
\DeclareMathOperator{\Ker}{Ker}
\DeclareMathOperator{\Spec}{Spec}

\DeclareMathOperator{\codim}{codim}
\DeclareMathOperator{\Img}{Im}

\begin{document}
\title{The Chow Motive of a Locally Trivial Fibration and Murre's conjectures}
\author{Carlos Pompeyo-Guti\'errez \\ carlos.pompeyo@ujat.mx \\ Divisi\'on Acad\'emica de Ciencias B\'asicas \\ Universidad Ju\'arez Aut\'onoma de Tabasco}
\date{February 2015}
%\classification{14C15}
%\keywords{Chow ring, Chow motives, Locally trivial fibration, Murre's conjectures.}
\maketitle
  
% Abstract
\begin{abstract}
Guillet and Soul\'e have shown in \cite{GIL96} that, for a fibration $\pi: Y \to X$ with fibre $Z$, locally trivial in the Zariski topology, we have a decomposition
\[ [Y] = [X] \cdot [Z], \]
where $[\cdot]$ denotes a class in the Grothendieck group $K_{0}(\cl{M}_{Rat}(k))$ associated to the category of (pure effective) Chow motives $\cl{M}_{Rat}(k)$ for a field $k$.
By assuming some additional properties for the fibre $Z$, we construct an explicit isomorphism $h(Y) \cong h(X) \otimes h(Z)$ in the category $\cl{M}_{Rat}(k)$, and we use it to prove, for this type of fibrations, some conjectures disscussed by Murre in  \cite{MUR93}.
\end{abstract}

% \maketitle is after abstract

\section{Introduction.}

In \cite{GIL96} Gillet and Soul\'e define, for any quasi-projective variety $X$ defined over a field $k$ of characteristic zero, a class $[X] \in K_{0}(\cl{M}_{Rat}(k))$ characterized by the following properties:
\begin{enumerate}
\item If $X$ is an smooth projective variety, then $[X]=[h(X)]$, where $h(X)=(X, \Delta_{X})$ denotes the Grothendieck motive under rational equivalence associated to $X$.
\item If $W \subset X$ is a closed subvariety in $X$,
\[ [X]=[W]+[X-W] \ . \]
\end{enumerate}

As a consequence, for a fibration $\pi: Y \to X$ with fiber $Z$, locally trivial for the Zariski topology of $X$, we have that
\[ [Y]=[X]\cdot[Z] \]
in $K_{0}(\cl{M}_{Rat}(k))$. Since $[X \times Z]=[X] \cdot [Z]$ in $K_{0}(\cl{M}_{Rat}(k))$, this implies that
\[ h(Y) \oplus M \cong \left( h(X) \otimes h(Z) \right) \oplus M \]
for some $M \in \ob(\cl{M}_{Rat}(k))$. Unfortunately, the cancellation law is not valid in general for the category $\cl{M}_{Rat}(k)$, as there are examples of fields $k$ for which the cancellation law fails, see for example, Remark 2.8 in \cite{CPSZ06}.

The main result of this work is the following:

\begin{theorem} \label{teo12}
Let $\pi : Y \to X$ be a locally trivial fibration with $\pi$ being a proper morphism and with fibres isomorphic to a fixed variety $Z$ having a Chow stratification. Then we have
\[ h(Y) \cong h(X) \otimes h(Z). \quad \finf \]
\end{theorem}

In Theorem \ref{teo12} we give an explicit isomorphism, closely related to the geometry of the fibre. After that, we discuss some conjectures proposed by Murre in \cite{MUR93}, namely we prove:

\begin{theorem} \label{thmCK}
Let $\pi:Y \to X$ a locally trivial fibration as in Theorem \ref{teo12}. Furthermore, suppose $X$ has a Chow-K\"unneth decomposition $\pi_{0}(X), \dots,$ $\pi_{2 \dim(X)}(X)$ satisfying that 
\[ \pi_{0}(X), \dots, \pi_{j-1}(X) , \pi_{2j+1}(X) , \dots , \pi_{2 \dim(X)}(X) \]
act as zero on $CH^{j}(X)$ ($0 \leq j \leq 2\dim(X)$). Then $Y$ has a Chow-K\"unneth decomposition $\pi_{0}(Y), \dots, \pi_{2 \dim(Y)}(Y)$ and for each $j$, 
\[ \pi_{0}(Y), \dots, \pi_{j-1}(Y), \pi_{2j+1}(Y), \dots, \pi_{2 \dim(Y)}(Y) \]
act as zero on $CH^{j}(Y)$. \fin
\end{theorem}

As an example of the previously mentioned fibrations we have the flag bundles associated to a given vector bundle $E$ over a variety $X$, more generally, any locally trivial fibration with fibres isomorphic to a quotient of a linear algebraic group by a parabolic subgroup. Our results are true even for varieties defined over a field having positive characteristic or not being algebraically closed.

\section{The Chow ring of a locally trivial fibration.}

In this section we calculate the Chow ring of certain locally trivial fibrations. The results given here are based in Proposition 1 in \cite{EDI97} and Lemma 2.8 in \cite{ELL89} and we include them for the sake of completeness.
We will consider fibrations with fibres isomorphic to a given variety $Z$ satisfying the following properties. 

\begin{definition} \label{chowpair}
Let $Z$ be a smooth projective variety with dimension $n$. We say that $Z$ satisfies the \textbf{Chow pairing conditions} if for each $p$ such that $0 \le p \le n$ we can find cycle classes $\tau_{p,1},...,\tau_{p,m_{p}} \in CH^{p}(Z)$ such that
\begin{enumerate}
\item $CH^{n}(Z) \cong \db{Z} \tau_{n,1}$ (where $\tau_{n,1}$ denotes the class of a point in $Z$).
\item For $p < n$, $CH^{p}(Z)$ is a free $\db{Z}$-module with finite rank
\[ CH^{p}(Z) \cong \bigoplus_{i=1}^{m_{p}} \db{Z} \tau_{p,i}. \]
\item For each $p < n$, we can give a perfect pairing
\[ CH^{p}(Z) \times CH^{n-p}(Z) \to CH^{n}(Z) \cong \db{Z}\tau_{n,1} \]
satisfying 
\[ \tau_{p,i} \cap \tau_{n-p,j} = \left\lbrace \begin{array}{ccc} \tau_{n,1} & if & i=j \\ 0 & if & i \neq j. \end{array} \right. \]
\end{enumerate}
\end{definition}

\begin{remark} Observe that both $CH^{p}(Z)$ and $CH^{n-p}(Z)$ have the same rank, in particular, $m_{0}=m_{n}=1$.
\end{remark}

\begin{definition} \label{chowstrat}
Let $Z$ be a smooth projective variety. We say that $Z$ has a \textbf{Chow stratification} if
\begin{enumerate}
 \item $Z$ has a cellular decomposition
\[ Z=Z_{d} \supset Z_{d-1} \supset \cdots \supset Z_{0} \supset Z_{-1}= \emptyset \]
by closed subvarieties such that each $Z_{i}-Z_{i-1}$ is a disjoint union of schemes $U_{i,j}$ isomorphic to affine spaces $\db{A}^{d_{i,j}}$.
\item $Z$ satisfies the Chow pairing conditions by taking the cycles appearing in the previous definition as $\tau_{i,j}:=\overline{U}_{i,j}$.
\end{enumerate}
\end{definition}

\begin{remark}. (See Example 1.9.1 in \cite[p. 23]{FUL98}). If a variety $Z$ has a cellular decomposition then $CH^{\ast}(Z)$ is finitely generated as a $\db{Z}$-module by the cycle classes $\overline{U}_{i,j}$; therefore the variety $Z$ will have a Chow stratification if it satisfies conditions $(i)$ and $(iii)$ in the definition of the Chow pairing conditions.
\end{remark}

It will prove to be useful to establish the following convention. 

\textbf{Convention.} For the rest of the paper, each time we say that $\pi: Y \to X$ is a locally trivial fibration, we will be assuming that $\pi$ is a locally trivial fibration with $\pi$ being a proper, flat, smooth morphism and with fibres isomorphic to a fixed variety $Z$ satisfying the Chow pairing conditions, unless otherwise stated.

Now, consider a locally trivial fibration $\pi:Y \to X$ and let $U \subset X$ be an open subset for which $\pi$ becomes trivial, and set $W:=X \setminus U$. Let $i:U \to X$, $j:W \to X$, $\imath : Y|_{U} \to Y$ and $\jmath : Y|_{W} \to Y$ be the inclusions, and denote by $\pi_{U}$ (resp. $\pi_{W}$) the restriction of $\pi$ to $Y|_{U}$ (resp. $Y|_{W}$). Let $\eta:Y|_{U} \to Z$ be the morphism induced by the projection of $U \times Z$ on the second factor.

Then for each $p$ we have the following diagram
\begin{align} \label{dia0301}
\begin{array}{c}
\xymatrix{ & & CH^{p}(Z) \ar@{>}[d]^{\eta^{\ast}} & \\
CH^{p-m}(Y|_{W}) \ar@{>}[r]^{\jmath_{\ast}} \ar@{>}[d]^{\pi_{W \ast}} & CH^{p}(Y) \ar@{>}[r]^{\imath^{\ast}} \ar@{>}[d]^{\pi_{\ast}} & CH^{p}(Y|_{U}) \ar@{>}[r] \ar@{>}[d]^{\pi_{U \ast}} & 0  \\
CH^{p-n-m}(W) \ar@{>}[r]_{j_{\ast}} & CH^{p-n}(X) \ar@{>}[r]_{i^{\ast}} & CH^{p-n}(U) \ar@{>}[r] & 0 } 
\end{array}
\end{align}
where $m$ denotes the codimension of $W$ in $X$; the rows are exact by Proposition $1.8$ in \cite{FUL98}, the left square commute by the functoriality of the push-forward and the right square commute by Proposition $1.7$ in \cite{FUL98}.

Using this diagram define elements $T_{p,i} \in CH^{p}(Y)$ such that $\imath^{\ast} T_{p,i} = \eta^{\ast} \tau_{p,i}$. We have the following Theorem.

\begin{theorem} (Duality Theorem). \label{teo0301}
Let $\pi:Y \to X$ be a fibration with fibre $Z$ satisfying the Chow pairing conditions. Then for any $p$, $q$ satisfying $p+q \leq n$ and any $\alpha \in CH^{\ast}(X)$ we have:
\[ \pi_{\ast}(\pi^{\ast}(\alpha) \cap T_{p,i} \cap T_{q,j}) = \left\lbrace \begin{array}{cl} \alpha & if \  (q,j)=(n-p,i) \\ 0 & otherwise \end{array} \right. \]
\end{theorem}

\proof By the projection formula
\[ \pi_{\ast}(\pi^{\ast}(\alpha) \cap T_{p,i} \cap T_{q,j})= \alpha \cap \pi_{\ast}(T_{p,i} \cap T_{q,j}),  \]
so it will be enough to calculate $\pi_{\ast}(T_{p,i} \cap T_{q,j})$. Observe that 
\[ \pi_{\ast}(T_{p,i} \cap T_{q,j}) \in CH^{p+q-n}(X), \]
and therefore $\pi_{\ast}(T_{p,i} \cap T_{q,j})=0$ if $p+q<n$. From now on we will suppose $q=n-p$. In this case, by looking at (\ref{dia0301}) we see that $CH^{0}(X) \cong CH^{0}(U)$.

Since the right square in (\ref{dia0301}) commutes then
\[ i^{\ast}\pi_{\ast}(T_{p,i} \cap T_{q,j})= \pi_{U \ast} \imath^{\ast}(T_{p,i} \cap T_{q,j}) \ . \]

But then
\[ \imath^{\ast}(T_{p,i} \cap T_{q,j})=\imath^{\ast}(T_{p,i}) \cap \imath^{\ast}(T_{q,j})=\eta^{\ast}(\tau_{p,i}) \cap \eta^{\ast}(\tau_{q,j}) = \eta^{\ast} (\tau_{p,i} \cap \tau_{q,j}) \ . \]

Now, if $j=i$ then $\tau_{p,i} \cap \tau_{q,j}=e$ and then 
\[ \imath^{\ast}(T_{p,i} \cap T_{q,j})= \eta^{\ast}(e)=1_{U} \times e , \]
where $e$ is the class of a point in $Z$ and therefore we have
\[ i^{\ast}\pi_{\ast}(T_{p,i} \cap T_{q,j})=\pi_{U \ast}(1_{U} \times e)=1_{U} \ . \]

So, being $i^{\ast}$ injective for $CH^{0}(X)$, we have $\pi_{\ast}(T_{p,i} \cap T_{q,j})=1_{X}$ if $j=i$, and therefore in this case
\[ \pi_{\ast}(\pi^{\ast}(\alpha) \cap T_{p,i} \cap T_{q,j})= \alpha \ . \]

Now suppose $j \neq i$, so we have $\tau_{p,i} \cap \tau_{q,j}=0$ and then
\[ \imath^{\ast}(T_{p,i} \cap T_{q,j})= \eta^{\ast}(0)=0 , \]
consequently,
\[ i^{\ast}\pi_{\ast}(T_{p,i} \cap T_{q,j})= \pi_{U \ast} \imath^{\ast}(T_{p,i} \cap T_{q,j})=0 \]
and since $i^{\ast}$ is injective for $CH^{0}(X)$
\[ \pi_{\ast}(T_{p,i} \cap T_{q,j})=0 \]
from which we conclude
\[ \pi_{\ast}(\pi^{\ast}(\alpha) \cap T_{p,i} \cap T_{q,j})=0 \]
for $j\neq i$. \fin

\textbf{Convention.} Since any element $\alpha_{i,j} \otimes n_{i,j} \tau_{i,j} \in \dis  CH^{p-i}(X) \otimes \db{Z} \tau_{i,j}$ can be rewritten as
\[ \alpha_{i,j} \otimes n_{i,j} \tau_{i,j} = n_{i,j}\alpha_{i,j} \otimes \tau_{i,j} = \beta_{i,j} \otimes \tau_{i,j} \ , \]
from now on we will simply denote such an element as $\alpha_{i,j} \otimes \tau_{i,j}$. \fin 

As a consequence of the Duality Theorem we have the following.
\begin{corollary} \label{cor0301}
 Let $\pi: Y \to X$ be a locally trivial fibration as before. Then
\[ \begin{array}{rrcl}
\fhi: & \dis \bigoplus_{i=0}^{p} \bigoplus_{j=1}^{m_{i}} \left( CH^{p-i}(X) \otimes \db{Z} \tau_{i,j} \right) & \to & CH^{p}(Y)  \\
 & \dis (\alpha_{i,j} \otimes \tau_{i,j})_{i,j}  & \mapsto & \dis \sum_{i=0}^{p} \sum_{j=1}^{m_{i}} \pi^{\ast}( \alpha_{i,j}) \cap T_{i,j} 
\end{array} \]
 is injective.
\end{corollary}

\proof Let $(\alpha_{i,j} \otimes \tau_{i,j})_{i,j} \in \ker \fhi$, so it satisfies the equation
\[ \sum \pi^{\ast}(\alpha_{i,j}) \cap T_{i,j} = 0. \]

Suppose we have $(\alpha_{i,j} \otimes \tau_{i,j})_{i,j} \neq 0$ and let $(k,l)$ be the (lexicographically)  greatest index such that $\alpha_{k,l} \otimes \tau_{k,l} \neq 0$. Multiplying the last equality by $T_{n-k,l}$ and then applying $\pi_{\ast}$ we obtain
\[ 0= \pi_{\ast} \left( \sum_{i,j} \pi^{\ast}(\alpha_{i,j}) \cap T_{i,j} \cap T_{n-k,l} \right) = \sum_{i,j} \pi_{\ast} (\pi^{\ast}(\alpha_{i,j}) \cap T_{i,j} \cap T_{n-k,l})=\alpha_{k,l} \]
which is absurd. Therefore $\ker \fhi =0$ and $\fhi$ is injective. \fin 

Is worth noticing that the group:
\[ \bigoplus_{i=0}^{p} \bigoplus_{j=1}^{m_{i}} \left( CH^{p-i}(X) \otimes \db{Z} \tau_{i,j} \right) \]
is isomorphic to the $p$-graded part of the graded ring $CH^{\ast}(X) \otimes CH^{\ast}(Z)$. In this way, Corollary \ref{cor0301} can be restated as follows.

\begin{corollary} \label{cor0302}
Let $\pi : Y \to X$ be a fibration as in Corollary \ref{cor0301}. Then $CH^{\ast}(X) \otimes CH^{\ast}(Z)$ is a $CH^{\ast}(X)$-submodule of $CH^{\ast}(Y)$. \fin
\end{corollary}

Now we center our attention on deciding when the morphism defined in Corollary \ref{cor0301} is surjective. In order to answer this question we require the following.

\begin{lemma} \label{lem0301}
If $Z$ has a Chow stratification, then for any variety $X$ we have
\[ CH^{\ast}(X \times Z) \cong CH^{\ast}(X) \otimes CH^{\ast}(Z). \]
\end{lemma}

\proof Notice that, in this case, the morphism from Corollary \ref{cor0301} can be written as
\[ \begin{array}{rrcl}
\fhi: & \dis \bigoplus_{i=0}^{p} CH^{p-i}(X) \otimes CH^{i}(Z) & \to & CH^{p}(X \times Z)  \\
 & \dis (\alpha_{i} \otimes \beta_{i})_{i}  & \mapsto & \dis \sum_{i=0}^{p} \pi_{X}^{\ast}(\alpha_{i}) \cap \pi_{Z}^{\ast}(\beta_{i})
\end{array} \]
where $\pi_{X}, \ \pi_{Z}$ are the projections from $X \times Z$ to $X$ and $Z$ respectively, and therefore this morphism is injective. Moreover we have the equality 
\[ \sum_{i} \pi_{X}^{\ast}(\alpha_{i}) \cap \pi_{Z}^{\ast}(\beta_{i}) = \sum_{i} \alpha_{i} \times \beta_{i} \ . \]
So, by comparing with Example 1.10.2 from \cite[p. 25]{FUL98} we obtain the surjectivity.
\fin

To conclude this Section, we provide the following Theorem.

\begin{theorem} \label{teo11}
Let $\pi : Y \to X$ be a locally trivial fibration with $\pi$ being a proper morphism and with fibres isomorphic to a fixed variety $Z$ having a Chow stratification. Then we have a $CH^{\ast}(X)$-modules isomorphism
\[ CH^{\ast}(Y) \cong CH^{\ast}(X) \otimes CH^{\ast}(Z). \quad \finf \]
\end{theorem}

\textbf{Proof.} At this point we only have to show that the morphism defined in Corollary \ref{cor0301} is surjective. In order to do this, we proceed by induction on the dimension of the base space $X$.

For $\dim X=0$ the result is trivial.

Now, suppose $\dim X >0$, let $U \subset X$ be an open set such that $Y$ becomes trivial on $U$, and set $W=X \setminus U$, $m=\codim_{X}(W)$. We have a diagram
\begin{align} \label{dia0302}
\begin{array}{c}
 \xymatrix{ \dis \bigoplus_{i=0}^{p-m} \bigoplus_{j=1}^{m_{i}} CH^{p-m-i}(W) \otimes \db{Z} \tau_{i,j} \ar@{>}[r]^<<<<{g'} \ar@{>}[d]_{h'} & CH^{p-m}(Y|_{W}) \ar@{>}[r] \ar@{>}[d]^{f'} & 0 \\ \dis \bigoplus_{i=0}^{p} \bigoplus_{j=1}^{m_{i}} CH^{p-i}(X) \otimes \db{Z} \tau_{i,j} \ar@{>}[r]^<<<<<<{\fhi} \ar@{>}[d]_{h} & CH^{p}(Y) \ar@{>}[d]^{f} & \\ \dis \bigoplus_{i=0}^{p} CH^{p-i}(U) \otimes CH^{i}(Z) \ar@{>}[r]^<<<<<<{g} \ar@{>}[d] & CH^{p}(U \times Z ) \ar@{>}[d] \ar@{>}[r] & 0 \\ 0 & 0 & }
 \end{array}
 \end{align}
where the first row is exact by induction hypothesis since $\dim W < \dim X$, the last row is given by Lemma \ref{lem0301} and the left column is obtained factor by factor by tensorizing the corresponding exact sequences obtained from Proposition $1.8$ in \cite{FUL98}.

Pick an element $\beta \in CH^{p}(Y)$ and set $\alpha_{1}:=f(\beta) $. Define elements $\alpha_{2}, \alpha_{3}$ satisfying $g(\alpha_{2})=\alpha_{1}$ and $h(\alpha_{3})=\alpha_{2}$. Then
\[ f(\fhi(\alpha_{3}))=g(h(\alpha_{3}))=g(\alpha_{2})=\alpha_{1}=f(\beta) \]
and therefore $\beta - \fhi(\alpha_{3}) \in \Ker f = \Img f'$, so we can write $\beta- \fhi(\alpha_{3})=f'(\alpha_{4})$ for some $\alpha_{4}$. Finally define $\alpha_{5}$ as an element satisfying  that $g'(\alpha_{5})=\alpha_{4}$. Then
\[ \fhi(\alpha_{3}+h'(\alpha_{5}))=\fhi(\alpha_{3})+f'(g'(\alpha_{5}))=\beta \]
and we have that $\fhi$ is surjective. \fin

\section{Defining the projectors.}

Along this section we will be using notations and conventions established by Manin in \cite{MAN68}. In this section we define pairwise orthogonal projectors for a locally trivial fibration as the one described in Section 2. In order to do this, we follow the ideas given by Manin in \cite{MAN68} to calculate the Chow motive of a projective bundle associated to a given vector bundle, and K\"ock in \cite{KCK91} which generalized the construction of Manin to  grassmannian bundles.

 We will define correspondences $p_{i,j} \in \ho_{CV(k)}(Y,Y)$ using the isomorphism described in Corollary \ref{cor0301}. In order to do this consider the auxiliary sets
\[ W_{i,j}= \{ (i,l) | j < l \leq m_{i} \} \cup \{ (k,l) | k>i, \ 1 \leq l \leq m_{k} \}, \]
and define the correspondences $p_{i,j}$ by a downward induction, starting with
\[ p_{n,m_{n}}=c_{T_{n,m_{n}}} \circ c(\pi) \circ c(\pi)^{t} \circ c_{T_{0,m_{n}}} \]
and in the general case by writing
\begin{equation} \label{ecp}
p_{i,j}=c_{T_{i,j}} \circ c(\pi) \circ c(\pi)^{t} \circ c_{T_{n-i,j}} \circ \left(         \Delta_{Y}-\sum_{(k,l) \in W_{i,j}} p_{k,l}  \right) \ . 
\end{equation}

The following Lemma will be used to prove some properties satisfied by the correspondences just defined.

\begin{lemma} \label{le0401}

$\dis (p_{i,j})_{e} \left( \sum_{r=0}^{p} \sum_{s=1}^{m_{r}} \pi^{\ast}(\alpha_{r,s}) \cap T_{r,s} \right)= \left\{ \begin{array}{ccl} \pi^{\ast}(\alpha_{i,j}) \cap T_{i,j} & if & i \leq p \\ 0 & if & i>p \ . \end{array}\right.$
\end{lemma}

\proof We will use a downward induction. First, observe that

$\begin{array}{cl}
 & \dis (p_{n,m_{n}})_{e} \left(\sum_{r=0}^{p} \sum_{s=1}^{m_{r}} \pi^{\ast}(\alpha_{r,s}) \cap T_{r,s} \right) \\
 
= & \dis m_{T_{n,m_{n}}} \left(\pi^{\ast} \left(\pi_{\ast} \left( \sum_{r=0}^{p} \sum_{s=1}^{m_{r}} \pi^{\ast}(\alpha_{r,s}) \cap T_{r,s} \cap T_{0,m_{n}} \right) \right) \right) \\

= & \dis \sum_{r=0}^{p} \sum_{s=1}^{m_{r}} m_{T_{n,m_{n}}} \left(\pi^{\ast} \left(\pi_{\ast} \left(  \pi^{\ast}(\alpha_{r,s}) \cap T_{r,s} \cap T_{0,m_{n}} \right) \right) \right)
\end{array}$

Now, by using Theorem \ref{teo0301}, last expression becomes
\[  m_{T_{n,m_{n}}} (\pi^{\ast}(\alpha_{n,m_{n}}) )= \pi^{\ast}(\alpha_{n,m_{n}})\cap T_{n,m_{n}}, \]
so we have verified the Lemma in this case.

Now, in order to clarify the general case, consider the sets
\[ M_{i,j}:=\{ (k,l) \ | \ k<i \} \cup \{ (i,l) \ | \ l \leq j \}. \]
Then, by applying the induction hypothesis:

$ \begin{array}{cl} &  \dis \left((\Delta_{Y})_{e}-\sum_{(k,l) \in W_{i,j}} (p_{k,l})_{e} \right) \left( \sum_{r=0}^{p} \sum_{s=1}^{m_{r}} \pi^{\ast}(\alpha_{r,s}) \cap T_{r,s} \right) \\ \ \\
= & \dis \sum_{r=0}^{p} \sum_{s=1}^{m_{r}} \pi^{\ast}(\alpha_{r,s}) \cap T_{r,s}-\sum_{(k,l) \in W_{i,j}} \pi^{\ast}(\alpha_{k,l}) \cap T_{k,l}  \\ \ \\
= & \dis \sum_{(r,s) \in M_{i,j}} \pi^{\ast}(\alpha_{r,s}) \cap T_{r,s} \\ \end{array} $

and so

$ \begin{array}{cl} &  \dis (p_{i,j})_{e} \left( \sum_{r=0}^{p} \sum_{s=1}^{m_{r}} \pi^{\ast}(\alpha_{r,s}) \cap T_{r,s} \right) \\
\ \\
= & \dis m_{T_{i,j}} \left( \pi^{\ast} \left( \pi_{\ast} \left(  \sum_{(r,s) \in M_{i,j}} \pi^{\ast}(\alpha_{r,s}) \cap T_{r,s} \cap T_{n-i,j} \right) \right) \right) \\
\ \\
= & \dis \sum_{(r,s) \in M_{i,j}} m_{T_{i,j}} \left( \pi^{\ast} \left( \pi_{\ast} \left(  \pi^{\ast}(\alpha_{r,s}) \cap T_{r,s} \cap T_{n-i,j} \right) \right) \right) . \\ \end{array} $

If $(r,s) \in M_{i,j}$ then $r \leq i$, and therefore $r+n-i \leq n$; then applying Theorem \ref{teo0301} we obtain
\[ \pi_{\ast} \left(  \pi^{\ast}(\alpha_{r,s}) \cap T_{r,s} \cap T_{n-i,j} \right)=\left\lbrace \begin{array}{cl} \alpha_{i,j} & if \  (r,s)=(i,j) \\ 0 & otherwise \end{array} \right. \]

In this way we can write
\[ \dis (p_{i,j})_{e} \left( \sum_{r=0}^{p} \sum_{s=1}^{m_{r}} \pi^{\ast}(\alpha_{r,s}) \cap T_{r,s} \right)=m_{T_{i,j}}(\pi^{\ast}(\alpha_{i,j}))=\pi^{\ast}(\alpha_{i,j}) \cap T_{i,j} \]
as desired. \fin

We will need the following two lemmas, the proof of which can be found in \cite{MAN68}.

\begin{lemma} \label{funcident}
For any morphism of varieties $\fhi:X \to Y$, any $T \in \ob(V(k))$ and any element $\alpha \in CH^{\ast}(X)$ we have
\begin{description}
\item[a)] $c(\fhi)_{T}=(id_{T} \times \fhi)^{\ast}$,
\item[b)] $c(\fhi)_{T}^{t}=(id_{T} \times \fhi)_{\ast}$,
\item[c)] $(c_{\alpha})_{T}=m_{1_{T} \times \alpha}$. \fin
\end{description}
\end{lemma}

\begin{lemma} \label{Yonconseq}
Let $\mathcal{D}$ be a diagram of objects and morphisms from the category $CV(k)$. Furthermore, let $I$ be
\[ I = \sum_{i=1}^{r} a_{i} f_{i}, \]
where $a_{i} \in \db{Z}$ and $f_{i}$ are some correspondences between the objects of the diagram $\mathcal{D}$. For $T \in \ob(V(k))$, let $I_{T}$ be
\[ I_{T} = \sum_{i=1}^{r} a_{i} (f_{i})_{T}. \]
Then $I=0$ if and only if $I_{T}=0$ for all $T \in \ob(V(k))$.  \fin
\end{lemma}

An immediate consequence of Lemmas \ref{funcident} and \ref{Yonconseq} is Manin's Identity Principle, which we state in what follows.

Suppose we have a diagram $\mathcal{D}$ of objects and morphisms of the category $V(k)$, and let $J$ be
\[ J = \sum_{i=1}^{r} a_{i}F_{i}, \]
where $a_{i} \in \db{Z}$ and every homomorphism $F_{i}$ is a composition of a finite number of homomorphisms  of the form $\fhi^{\ast}$, $\fhi_{\ast}$, $m_{\alpha}$ for $\alpha \in C(X)$, $X \in \ob(\mathcal{D})$, $\fhi \in \mor(\mathcal{D})$.

For any $T \in \ob(V(k))$ we denote by $T \times J$ the identity obtained from $J$ by changing all the objects $X$ by $T \times X$, all the morphisms $\fhi$ by $id_{T} \times \fhi$ and all the morphisms $m_{\alpha}$ by $m_{1_{T} \times \alpha}$.

In a similar way, denote by $c(J)$ the identity obtained from $J$ by changing all the morphisms $\fhi^{\ast}$ by $c(\fhi)$, all the morphisms $\fhi_{\ast}$ by $c(\fhi)^t$ and all the morphisms $m_{\alpha}$ by $c_{\alpha}$.

\begin{theorem} \label{Manin} \textbf{Manin's Identity Principle} (\cite[p.450]{MAN68}).
Let $J$ be as before. The following two assertions are equivalent.
\begin{description}
\item[a)] $T \times J=0$ for all $T \in \ob(V(k))$.
\item[b)] $c(J)=0$. \fin
\end{description}
\end{theorem}

The correspondences $p_{i,j}$ have the following properties.

\begin{theorem} \label{teo0401}
Let $p_{i,j}$ be the correspondences defined before. Then we have the following:
\begin{enumerate}
\item The correspondences $p_{i,j}$ are of degree zero.
\item $ \dis \sum_{i,j} p_{i,j} = \Delta_{Y} $
\item $ \dis p_{i,j} \circ p_{k,l}= \delta_{(i,j)}^{(k,l)} p_{i,j} $
\end{enumerate}
\end{theorem}

\proof The first affirmation is clear from the definition of the correspondences $p_{i,j}$. In order to prove the remaining assertions we will use Manin's Identity Principle. 

We are going to make more precise what we are supposed to prove. Define morphisms $\rho_{i,j}:CH^{\ast}(Y) \to CH^{\ast}(Y)$ in a inductive way by:
\[ \rho_{i,j}:=m_{T_{i,j}} \circ \pi^{\ast} \circ \pi_{\ast} \circ m_{T_{n-i,j}} \circ \left( id_{CH^{\ast}(Y)} - \sum_{(k,l) \in W_{i,j}} \rho_{k,l} \right),  \]
and for a variety $T$, denote by $\rho_{i,j}^{T}$ the morphism:
\[ \rho_{i,j}^{T}:=m_{1_{T} \times T_{i,j}} \circ (id_{T} \times \pi)^{\ast} \circ (id_{T} \times \pi)_{\ast} \circ m_{1_{T} \times T_{n-i,j}} \circ \left( id_{CH^{\ast}(T \times Y)} - \sum_{(k,l) \in W_{i,j}} \rho_{k,l}^{T} \right) \]

 Manin's Identity Principle asserts that identities $2$ and $3$ of this Theorem hold if and only if the identities
\begin{equation} \label{ec0403}
 \sum_{i,j} \rho_{i,j}^{T} = id_{CH^{\ast}(T \times Y)} 
\end{equation}
 and 
\begin{equation} \label{ec0404}
 \rho_{i,j}^{T} \circ \rho_{k,l}^{T}= \delta_{(i,j)}^{(k,l)} \rho_{i,j}^{T} 
\end{equation}
 hold for every variety $T$.
 
 Now  if
 \[ p_{i,j}^{T}=c_{1_{T} \times T_{i,j}} \circ c(id_{T} \times \pi) \circ c(id_{T} \times \pi)^{t} \circ c_{1_{T} \times T_{n-i,j}} \circ \left( \Delta_{T \times Y}-\sum_{(k,l) \in W_{i,j}} p_{k,l}^{T}  \right) \]
 then $(p_{i,j}^{T})_{e}=\rho_{i,j}^{T}$ and we can rewrite (\ref{ec0403}) and (\ref{ec0404}) as
 \begin{equation} \label{ec0405}
 \begin{array}{c}
  \dis \sum_{i,j} (p_{i,j}^{T})_{e} = (\Delta_{T \times Y})_{e} \ , \\ 
  \ \\
   (p_{i,j}^{T})_{e} \circ (p_{k,l}^{T})_{e}= \delta_{(i,j)}^{(k,l)} (p_{i,j}^{T})_{e} 
   \end{array}
 \end{equation}
 
Now, by Lemma \ref{le0401}
 \[ \dis \sum_{i=0}^{n} \sum_{j=1}^{m_{i}} (p_{i,j})_{e}\left( \sum_{r=0}^{p} \sum_{s=1}^{m_{r}} \pi^{\ast}(\alpha_{r,s}) \cap T_{r,s} \right) = \sum_{i=0}^{p} \sum_{j=1}^{m_{i}} \pi^{\ast}(\alpha_{i,j}) \cap T_{i,j} \]
so we have that $\dis \sum_{i,j}(p_{i,j})_{e}|_{CH^{p}(Y)}=id_{CH^{p}(Y)}$, and therefore
\begin{equation} \label{ec0406}
\sum_{i,j}(p_{i,j})_{e}=id_{CH^{\ast}(Y)}=(\Delta_{Y})_{e} \ . 
\end{equation}

On the other hand,
\[ (p_{i,j})_{e} \circ (p_{k,l})_{e} \left( \sum_{r=0}^{p} \sum_{s=1}^{m_{r}} \pi^{\ast}(\alpha_{r,s}) \cap T_{r,s} \right) = (p_{i,j})_{e} \left( \pi^{\ast}(\alpha_{k,l}) \cap T_{k,l} \right) \]
but
\[ (p_{i,j})_{e} \left( \pi^{\ast}(\alpha_{k,l}) \cap T_{k,l} \right) = \left\lbrace \begin{array}{ccl} \pi^{\ast}(\alpha_{i,j}) \cap T_{i,j} & if & (i,j)=(k,l) \\ 0 & if & (i,j) \ne (k,l) \end{array} \right. \]
therefore
\begin{equation} \label{ec0407}
   (p_{i,j})_{e} \circ (p_{k,l})_{e}= \delta_{(i,j)}^{(k,l)} (p_{i,j})_{e} \ .
\end{equation}

%and, by equation (\ref{ec0402})
%\[ \rho_{i,j} \circ \rho_{k,l}=(p_{i,j})_{e} \circ (p_{k,l})_{e}= \delta_{(i,j)}^{(k,l)} %(p_{i,j})_{e} = \delta_{(i,j)}^{(k,l)} \rho_{i,j}. \]

But both (\ref{ec0406}) and (\ref{ec0407}) are true for a locally trivial fibration $\pi:Y \to X$ satisfying the hypothesis of Theorem \ref{teo11}, and the locally trivial fibration 
\[ id_{T} \times \pi: T \times Y \to T \times X \]
satisfy such hypothesis, provided we can show that the elements $1_{T} \times T_{i,j}$ generate the Chow ring $CH^{\ast}(T \times Y)$ as a $CH^{\ast}(T \times X)$-module. But this follows from the definition of the mentioned generators if we do suitable changes to diagram (\ref{dia0301}). Therefore, identities (\ref{ec0405}) also hold. \fin 

Now we will establish some Lemmas needed to prove Theorem \ref{teo12}. From now on, each time we say we have an isomorphism of motives, we are refering to the fact that we have an isomorphism between the additive structures of the motives involved. We start by calculating some factors appearing in the decomposition that we will give later for the motive $h(Y)$.

\begin{lemma} \label{le0402}
Let $\pi : Y \to X$ and $\pi': Y' \to X$ be two locally trivial fibrations satisfying the hypothesis of Theorem \ref{teo11}. Then we have an isomorphism of motives
\[ h(Y) \cong h(Y'). \]
\end{lemma}
\proof Denote by $T_{i,j}$ (resp. $T_{i,j}'$) the generators of $CH^{\ast}(Y)$ (resp. $CH^{\ast}(Y')$) as a $CH^{\ast}(X)$-module.

Theorem \ref{teo0401} lets us conclude that
\[ h(Y)=(Y, \Delta_{Y}) = \left( Y, \sum_{i,j} p_{i,j} \right) \cong \bigoplus_{i,j} (Y,p_{i,j}) \]
and
\[ h(Y') \cong \bigoplus_{i,j} (Y',p_{i,j}'), \]
where $p_{i,j}, \ p_{i,j}'$ are defined as in (\ref{ecp}), so it will be enough to show that the factors appearing in these decompositions are isomorphic.

In order to do this, define morphisms of motives
\[h_{i,j} \in \ho_{\cl{M}_{Rat}(k)}((Y,p_{i,j}),(Y',p_{i,j}'))\]
by the formula:
\[ h_{i,j}:=c_{T_{i,j}'} \circ c(\pi') \circ c(\pi)^{t} \circ c_{T_{n-i,j}} \circ \left( \Delta_{Y} - \sum_{(k,l) \in W_{i,j}} p_{k,l} \right) \ ; \]
analogously define the morphisms $h_{i,j}' \in \ho_{\cl{M}_{Rat}(k)}((Y',p_{i,j}'),(Y,p_{i,j}))$ by:
\[ h_{i,j}':=c_{T_{i,j}} \circ c(\pi) \circ c(\pi')^{t} \circ c_{T_{n-i,j}'} \circ \left( \Delta_{Y'} - \sum_{(k,l) \in W_{i,j}} p_{k,l}' \right) \ . \]

At this point we would have to show the commutativity of the diagrams
\[ \xymatrix{ Y \ar@{>}[d]_{p_{i,j}} \ar@{>}[r]^{h_{i,j}} & Y' \ar@{>}[d]^{p_{i,j}'} & & Y' \ar@{>}[d]_{p_{i,j}'} \ar@{>}[r]^{h_{i,j}'} & Y \ar@{>}[d]^{p_{i,j}} \\ Y \ar@{>}[r]_{h_{i,j}} & Y' & & Y' \ar@{>}[r]_{h_{i,j}'} & Y } \]
but this follows by Manin's Identity Principle since we have both
\[ (h_{i,j})_{e} \left( \sum_{r=0}^{p} \sum_{s=1}^{m_{r}} \pi^{\ast}(\alpha_{r,s}) \cap T_{r,s} \right)= \pi'^{\ast}(\alpha_{i,j}) \cap T_{i,j}' \]
and a similar equation holding for $(h_{i,j}')_{e}$.

By following this procedure we can also obtain the identities
\[ h_{i,j}' \circ h_{i,j} = \Delta_{Y} \mod p_{i,j}, \quad h_{i,j} \circ h_{i,j}' = \Delta_{Y'} \mod p_{i,j}' \]
which expose both $h_{i,j}$ and $h'_{i,j}$ as the desired isomorphisms. \fin

\begin{lemma} \label{le0403}
Let $\pi: Y \to X$ be a locally trivial fibration satisfying the hypothesis of Theorem \ref{teo11}. Then we have an isomorphism of motives
\[ (Y, p_{i,j}) \cong h(X) \otimes ( Z, p_{i,j,Z}), \]
where the projectors $p_{i,j,Z} \in \ho_{\cl{M}_{Rat}(k)}(Z,Z)$ are defined by the formula
\[ p_{i,j,Z}:=\tau_{n-i,j} \times \tau_{i,j}. \]
\end{lemma}

\proof By Lemma \ref{le0402} we have that
\[ (Y, p_{i,j}) \cong (X \times Z, q_{i,j} ) \]
where the projectors $q_{i,j}$ are the ones defined for the trivial fibration
\[ \xymatrix{X \times Z \ar@{>}[r]^<<<<<{\rho} & X} \]
by using the formula (\ref{ecp}).

We will show that
\[ (X \times Z, q_{i,j} ) = ( X \times Z, \Delta_{X} \otimes p_{i,j,Z}) \]

We start by observing that the collection of correspondences $p_{i,j,Z}$ are, in fact, projectors in the category of motives. Clearly, the correspondences $p_{i,j,Z}$ have degree zero. Now, 
\begin{eqnarray}
 p_{i,j,Z} \circ p_{k,l,Z} & = & \pi_{13 \ast} \left(\pi_{12}^{\ast}(\tau_{n-i,j} \times \tau_{i,j}) \cap \pi_{23}^{\ast}(\tau_{n-k,l} \times \tau_{k,l}) \right) \nonumber \\
  & = & \pi_{13 \ast}(\tau_{n-i,j} \times (\tau_{i,j} \cap \tau_{n-k,l}) \times \tau_{k,l}) \ . \nonumber
\end{eqnarray}

Suppose $\tau_{i,j} \cap \tau_{n-k,l} \ne 0$. Then
\[ \pi_{13}(\tau_{n-i,j} \times (\tau_{i,j} \cap \tau_{n-k,l}) \times \tau_{k,l})=\tau_{n-i,j} \times \tau_{k,l} \ . \]

Since
\[ \dim(\tau_{n-i,j} \times (\tau_{i,j} \cap \tau_{n-k,l}) \times \tau_{k,l})=n,  \]
\[ \dim(\tau_{n-i,j} \times \tau_{k,l})=n+i-k  \]
we see that they have the same dimension if and only if $k=i$. Therefore
\[ \pi_{13 \ast}(\tau_{n-i,j} \times (\tau_{i,j} \cap \tau_{n-k,l}) \times \tau_{k,l})=\left \{ \begin{array}{cll} \tau_{n-i,j} \times \tau_{i,l} & if & k=i \\ 0 & if & k \ne i \end{array} \right. \]
but since we are assuming $\tau_{i,j} \cap \tau_{n-i,l} \ne 0$, we have that $l=j$, so
\[ \pi_{13 \ast}(\tau_{n-i,j} \times (\tau_{i,j} \cap \tau_{n-k,l}) \times \tau_{k,l})=\delta_{(i,j)}^{(k,l)} \tau_{n-i,j} \times \tau_{i,j} \]

and therefore
\[ p_{i,j,Z} \circ p_{k,l,Z}= \delta_{(i,j)}^{(k,l)} p_{i,j,Z} \ . \]

Now consider the case when $\tau_{i,j} \cap \tau_{n-k,l}=0$. Then $(i,j) \ne (k,l)$, otherwise we would have $\tau_{i,j} \cap \tau_{n-k,l}=e$; therefore $\delta_{(i,j)}^{(k,l)}=0$. Another consequence of $\tau_{i,j} \cap \tau_{n-k,l}=0$ is that, by Proposition $1.10$ in \cite[p. 24]{FUL98} we have that
\[ \tau_{n-i,j} \times (\tau_{i,j} \cap \tau_{n-k,l}) \times \tau_{k,l}=0. \]

So in this case we also have the equality
\[ p_{i,j,Z} \circ p_{k,l,Z}= \delta_{(i,j)}^{(k,l)} p_{i,j,Z} \ . \]

Therefore the correspondences $p_{i,j,Z}$ define mutually orthogonal projectors on $Z$. Now we proceed to verify that the projectors $q_{i,j}$ and $\Delta_{X} \otimes p_{i,j,Z}$ coincide.

By Lemma \ref{le0401} we have that
\[ (q_{i,j})_{e} \left( \sum_{r=0}^{p} \sum_{s=1}^{m_{r}} \rho^{\ast}(\alpha_{r,s}) \cap 1_{X} \times \tau_{r,s}  \right) = \rho^{\ast}(\alpha_{i,j}) \cap 1_{X} \times \tau_{i,j} \ . \]

On the other hand
\begin{eqnarray}
 (\Delta_{X} \otimes p_{i,j,Z})_{e}(\rho^{\ast}(\alpha_{r,s}) \cap 1_{X} \times \tau_{r,s}) & = & (p_{i,j,Z})_{X}(\rho^{\ast}(\alpha_{r,s}) \cap 1_{X} \times \tau_{r,s}) \nonumber \\
& = & p_{i,j,Z} \circ (\rho^{\ast}(\alpha_{r,s}) \cap 1_{X} \times \tau_{r,s}) \nonumber \\
& = & p_{13 \ast}( \alpha_{r,s} \times (\tau_{r,s} \cap \tau_{n-i,j}) \times \tau_{i,j}) \nonumber \\
& = & \left\lbrace \begin{array}{ccc} \alpha_{i,j} \times \tau_{i,j} & if & (r,s)=(i,j) \\ 0 & if & (r,s) \ne (i,j) \end{array} \right. \nonumber \\
& = & \delta_{(i,j)}^{(k,l)} \rho^{\ast}(\alpha_{i,j}) \cap 1_{X} \times \tau_{i,j} \ . \nonumber
\end{eqnarray}

Therefore
\begin{eqnarray}
& & (\Delta_{X} \otimes p_{i,j,Z})_{e} \left( \sum_{r=0}^{p} \sum_{s=1}^{m_{r}} \rho^{\ast}(\alpha_{r,s}) \cap 1_{X} \times \tau_{r,s}  \right) \nonumber \\
& = & \rho^{\ast}(\alpha_{i,j}) \cap 1_{X} \times \tau_{i,j} \ . \nonumber
\end{eqnarray}
At this point we have proved the equality
\[ (q_{i,j})_{e}=(\Delta_{X} \otimes p_{i,j,Z})_{e} \]
and by applying Manin's Identity Principle we obtain
\[ q_{i,j}=\Delta_{X} \otimes p_{i,j,Z} \ . \]

We conclude the proof of this Lemma by noticing that
\[ (Y,p_{i,j}) \cong (X \times Z, q_{i,j})=(X \times Z , \Delta_{X} \otimes p_{i,j,Z}) \cong (X, \Delta_{X}) \otimes (Z, p_{i,j,Z}) \]
as desired. \fin

\begin{lemma} \label{le0404}
Under the hypothesis of Lemma \ref{le0403} we have that
\[ h(Z) \cong \bigoplus_{i,j} (Z,p_{i,j,Z}) \]
\end{lemma}

\proof We have already shown that the correspondences $p_{i,j,Z}$ induce mutually orthogonal projectors. So, our proof will be finished if we can show that
\[ \sum_{i,j} p_{i,j,Z} = \Delta_{Z} \ . \]

The elements of $CH^{\ast}(Z)$ can be written as
\[ \sum_{r=0}^{n} \sum_{s=1}^{m_{r}} n_{r,s} \tau_{r,s} \]
for some $n_{r,s} \in \db{Z}$. Observe that
\[ (p_{i,j,Z})_{e}\left(\sum_{r=0}^{n} \sum_{s=1}^{m_{r}} n_{r,s} \tau_{r,s} \right)= p_{i,j,Z} \circ \left(\sum_{r=0}^{n} \sum_{s=1}^{m_{r}} n_{r,s} \tau_{r,s} \right)= \sum_{r=0}^{n} \sum_{s=1}^{m_{r}} n_{r,s} p_{i,j,Z} \circ \tau_{r,s}  \]
besides
\begin{eqnarray}
p_{i,j,Z} \circ \tau_{r,s} & = & p_{2 \ast}( (\tau_{n-i,j} \times \tau_{i,j}) \cap (\tau_{r,s} \times 1_{Z}) ) \nonumber \\
 & = & p_{2 \ast}( (\tau_{n-i,j} \cap \tau_{r,s}) \times \tau_{i,j} ) \nonumber \\
 & = & \left\lbrace \begin{array}{ccc}  \tau_{i,j} & if & (r,s)=(i,j) \\ 0 & if & (r,s) \ne (i,j) \end{array} \right. \nonumber \\
 & = & \delta_{(i,j)}^{(r,s)} \tau_{i,j} \nonumber
\end{eqnarray}
in this way we obtain that
\[(p_{i,j,Z})_{e}\left( \sum_{r=0}^{n} \sum_{s=1}^{m_{r}} n_{s} \tau_{r,s} \right)=n_{i,j} \tau_{i,j} \]
and therefore
\[ \sum_{i,j} (p_{i,j,Z})_{e}=id_{CH^{\ast}(Z)}=(\Delta_{Z})_{e} \]
and by using Manin's Identity Principle we obtain the desired result. \fin

Now we have at our disposal all the tools needed to prove Theorem \ref{teo12}. 

\textbf{Proof of Theorem \ref{teo12}} We have that
\[ \begin{array}{rcl} h(Y) = (Y, \Delta_{Y}) & \cong & \dis \bigoplus_{i,j} (Y, p_{i,j}) \\
\ \\
& \cong & \dis \bigoplus_{i,j} \left( h(X) \otimes (Z,p_{i,j,Z}) \right) \\
\ \\
& \cong & h(X) \otimes \left( \dis \bigoplus_{i,j} (Z,p_{i,j,Z}) \right) \\
\ \\
& \cong & h(X) \otimes h(Z). \end{array} \finf \]

\section{Murre's conjectures.}

We begin this section by recalling some definitions.

\begin{definition}
Let $X$ be an smooth $d$ dimensional projective variety over a field $k$. We say that $X$ has a Chow-K\"unneth decomposition if we can find cycle classes
\[ \pi_{0}(X), \dots , \pi_{2d}(X) \in CH^{d}(X \times X, \db{Q}) \] 
such that
\begin{description}
\item[a)] $\pi_{i}(X) \circ \pi_{j}(X) = \delta_{i,j} \pi_{i}(X)$.
\item[b)] $\dis \Delta_{X}= \sum_{i=0}^{2d} \pi_{i}(X)$.
\item[c)] (over $\bar{k}$) $\pi_{i}$ modulo (co)homological equivalence (for example, in \'etale cohomology) is the usual K\"unneth component $\Delta_{X}(2d-i,i)$.
\end{description}
If we define $h^{i}(X):=(X, \pi_{i}(X))$, then we will say that
\[ h(X)= \bigoplus_{i=0}^{2d} h^{i}(X) \]
(or equivalently, the collection $\pi_{0}(X), \dots, \pi_{2d}(X)$) is a Chow-K\"unneth (CK) decomposition for $X$.
\end{definition}

Some examples of varieties having a CK decomposition are curves \cite{KLE70}, surfaces \cite{MUR90}, products of curves and surfaces \cite{MUR93II}, abelian varieties \cite{SHE74} and uniruled complex 3-folds \cite{ANG98}. The following conjectures (among others) were proposed by Murre in \cite{MUR93}, and they are related to a conjectural filtration on the Chow groups of an algebraic variety.

\textbf{Murre's conjectures.} 
\begin{description}
\item[A)] Every smooth projective $d$ dimensional variety $X$ has a Chow-K\"unneth decomposition:
\[ h(X) \cong \bigoplus_{i=0}^{2d} (X, \pi_{i}(X)) \]
\item[B)] For each $j$, $\pi_{0}(X), \dots, \pi_{j-1}(X), \pi_{2j+1}(X), \dots, \pi_{2d}(X)$ act as zero on $CH^{j}(X, \db{Q})$.
\end{description}

In order to say something about the conjectures in case of the fibrations studied in this work, we need to establish some identities, the proof is straightforward.

\begin{lemma}
Let $X, Y, Z \in \ob(V(k))$, $\alpha \in CH^{\ast}(Y)$, $\fhi \in CH^{\ast}(X \times Y)$, $\psi \in CH^{\ast}(Y \times Z)$, $\tau \in CH^{\ast}(X \times Z)$. Let $f:X \to Y$ and $g:Y \to Z$ be morphisms in $V(k)$. Then we have the following identities.
\begin{enumerate}
\item $c_{\alpha} \circ \fhi = (1_{X} \times \alpha) \cdot \fhi$,
\item $\psi \circ c_{\alpha}= (\alpha \times 1_{X}) \cdot \psi$,
\item $c(f) \circ \fhi = (id_{X} \times f)^{\ast}(\fhi)$,
\item $c(g)^{t} \circ \fhi = (id_{X} \times g)_{\ast}(\fhi)$,
\item $\tau \circ c(f) = (f \times id_{Z})_{\ast}(\tau)$,
\item $\psi \circ c(f)^{t} = (f \times id_{Z})^{\ast}(\psi)$. \fin
\end{enumerate}
\end{lemma}

In particular, we can rewrite some parts of the correspondences $p_{i,j}$ given in the last section, namely
\[ c_{T_{i,j}} \circ c(\pi) \circ c(\pi)^{t} \circ c_{T_{n-i,j}}= (\pi \times \pi)^{\ast}(\Delta_{X}) \cdot (T_{n-i,j} \times T_{i,j}) \ . \]
From this new expression we see that, provided we have a CK decomposition for the base space of the fibration as
\[ \Delta_{X} = \sum_{i=0}^{2 \dim(X)} \pi_{i}(X) \ , \]
then a CK decomposition for the fibration $Y$ must involve correspondences with terms of the form
\[ (\pi \times \pi)^{\ast}(\pi_{r}(X)) \cdot (T_{n-i,j} \times T_{i,j})=c_{T_{i,j}} \circ c(\pi) \circ \pi_{r}(X) \circ c(\pi)^{t} \circ c_{T_{n-i,j}} \ . \]

That will be the case, as we will see in the proof of Theorem \ref{thmCK}. 

\textbf{Proof of Theorem \ref{thmCK}}. Let $\pi_{0}(X), \cdots, \pi_{2d}(X)$ be a CK decomposition for $X$, where $d$ is the dimension of $X$. For a cycle $\fhi \in CH^{\ast}(X \times X)$, define a cycle $\fr{p}_{j}(\fhi) \in CH^{\ast}(Y \times Y)$ as (here $n=\dim(Z)$)

\[ \dis \fr{p}_{j}(\fhi)= \sum_{\lambda=1}^{m_{j/2}} c_{T_{j/2, \lambda}} \circ c(\pi) \circ \fhi \circ c(\pi)^{t} \circ c_{T_{n-j/2, \lambda}} \circ \left( \Delta_{Y} - \sum_{(k,l) \in W_{j/2, \lambda}} p_{k,l} \right) \]
for $j$ even and $\fr{p}_{j}(\fhi)=0$ for $j$ odd.

Since
\[ c_{T_{j/2, \lambda}} \circ c(\pi) \circ \fhi \circ c(\pi)^{t} \circ c_{T_{n-j/2, \lambda}} = (\pi \times \pi)^{\ast}(\fhi) \cdot (T_{n-j/2, \lambda} \times T_{j/2, \lambda}) \ , \]
it follows that $\fr{p}_{j}(\cdot)$ is additive.

Now, for each integer $k$ between $0$ and $2 \dim(Y)$ define the set
\[ I_{k}:=\{ (i,j) \ | \ 0 \leq i \leq 2d, \ 0 \leq j \leq 2n, \ i+j=k \} \]
and the correspondence
\[ \pi_{k}(Y):= \sum_{(i,j) \in I_{k}} \fr{p}_{j}(\pi_{i}(X)). \]

We will show that $\pi_{0}(Y), \dots, \pi_{2 \dim(Y)}(Y)$ give a CK decomposition for $Y$ and satisfy the conjecture \textbf{B)} established before.

Clearly $\pi_{k}(Y)$ is a correspondence of degree zero. Besides, since $I_{k} \cap I_{k'} = \emptyset$ for $k \ne k'$ and 
\[ \bigcup_{k=0}^{2 \dim(Y)} I_{k} = \{0, 1, \dots, 2d \} \times \{0,1, \dots, 2n \} \] 
it follows that
\begin{eqnarray}
\sum_{k=0}^{2 \dim(Y)} \pi_{k}(Y) & = & \sum_{k=0}^{2 \dim(Y)} \sum_{(i,j) \in I_{k}} \fr{p}_{j}(\pi_{i}(X)) \nonumber \\ & = & \sum_{i=0}^{2d} \sum_{j=0}^{2n} \fr{p}_{j}(\pi_{i}(X)) \nonumber \\ & = & \sum_{j=0}^{2n} \fr{p}_{j} \left( \sum_{i=0}^{2d} \pi_{i}(X) \right) \nonumber \\ & = & \sum_{j=0}^{n} \sum_{\lambda=0}^{m_{j}} p_{j, \lambda} = \Delta_{Y} \ . \nonumber
\end{eqnarray}

Now, we will verify that we have the identities
\[ \pi_{k}(Y) \circ \pi_{k'}(Y) = \delta_{k,k'} \pi_{k}(Y) \ . \]

By Lemma \ref{Yonconseq}, we have to show that for any variety $T \in \ob(V(k))$,
\[ \pi_{k}(Y)_{T} \circ \pi_{k'}(Y)_{T} = \delta_{k,k'} \pi_{k}(Y)_{T} \ . \]

As we observed in the proof of Theorem \ref{teo0301}, $id_{T} \times \pi: T \times Y \to T \times X$ is a locally trivial fibration with fibres isomorphic to $Z$ and having a Chow stratification, the generators of $CH^{\ast}(T \times Y)$ as $CH^{\ast}(T \times X)$-module being the elements $1_{T} \times T_{i,j}$. Therefore, by Theorem \ref{teo11} any element $\beta \in CH^{p}(T \times Y)$ can be written as
\[ \beta = \sum_{r=0}^{p} \sum_{s=1}^{m_{r}} (id_{T} \times \pi)^{\ast}(\alpha_{r,s}) \cap (1_{T} \times T_{r,s}) \ , \]
for some $\alpha_{r,s} \in CH^{p-r}(T \times X)$.

By Lemma \ref{funcident}, each non zero summand of $\fr{p}_{j}(\pi_{i}(X))_{T}$ is of the form
\[ m_{1_{T} \times T_{j/2, \lambda}} \circ (id_{T} \times \pi)^{\ast} \circ \pi_{i}(X)_{T} \circ (id_{T} \times \pi)_{\ast} \circ m_{1_{T} \times T_{n-j/2, \lambda}} \circ \left( id_{T \times Y} - \sum_{(k,l) \in W_{j/2, \lambda}} (p_{k,l})_{T} \right)  \]
and by doing similar calculations to the ones given in the proof of Lemma \ref{le0401} , we obtain
\[ \left( (id_{T} \times \pi)_{\ast} \circ m_{1_{T} \times T_{n-j/2, \lambda}} \circ \left( id_{T \times Y} - \sum_{(k,l) \in W_{j/2, \lambda}} (p_{k,l})_{T} \right) \right)(\beta)= \alpha_{j/2, \lambda} \ ; \]
in consequence
\[ \fr{p}_{j}(\pi_{i}(X))_{T}(\beta)= \dis \sum_{\lambda=1}^{m_{j/2}} (id_{T} \times \pi)^{\ast}(\pi_{i}(X)_{T}(\alpha_{j/2, \lambda})) \cap (1_{T} \times T_{j/2, \lambda}) \ ; \]
observe that in this last expression, by applying $\fr{p}_{j}(\pi_{i}(X))_{T}$ to an element, the result only involve terms in which the generator is of the form $1_{T} \times T_{j/2, \lambda}$. Therefore, should we apply $\fr{p}_{j'}(\pi_{i'}(X))_{T}$ to $\fr{p}_{j}(\pi_{i}(X))_{T}(\beta)$ for $j' \ne j$ (and no matter what value of $i'$ we choose) we would obtain zero. On the other hand, if $j=j'$ then
\[ \fr{p}_{j}(\pi_{i'}(X))_{T} \left(\fr{p}_{j}(\pi_{i}(X))_{T}(\beta) \right)= \dis \sum_{\lambda=1}^{m_{j/2}} (id_{T} \times \pi)^{\ast}( (\pi_{i'}(X) \circ \pi_{i}(X))_{T}(\alpha_{j/2, \lambda})) \cap (1_{T} \times T_{j/2, \lambda}) \ ; \]
but we have that $\pi_{i'}(X) \circ \pi_{i}(X) = \delta_{i',i} \pi_{i}(X)$, so it follows that
\[ \fr{p}_{j}(\pi_{i'}(X))_{T} \left(\fr{p}_{j}(\pi_{i}(X))_{T}(\beta) \right)= \delta_{i',i} \fr{p}_{j}(\pi_{i}(X))_{T}(\beta) \ . \]

To summarize, we have that
\[ \fr{p}_{j'}(\pi_{i'}(X))_{T} \left(\fr{p}_{j}(\pi_{i}(X))_{T}(\beta) \right) = \delta_{(i,j)}^{(i',j')} \fr{p}_{j}(\pi_{i}(X))_{T}(\beta) \ . \]

Now, when $k \ne k'$, $I_{k} \cap I_{k'} = \emptyset$, and therefore $\delta_{(i,j)}^{(i',j')}=0$ for any $(i,j) \in I_{k}$ and any $(i',j') \in I_{k'}$. Therefore, for $k \ne k'$
\[ \pi_{k'}(Y)_{T} \circ \pi_{k}(Y)_{T}(\beta)= 0 \ . \]

In a similar fashion, $\pi_{k}(Y)_{T} \circ \pi_{k}(Y)_{T}(\beta)= \pi_{k}(Y)_{T}(\beta) $. Therefore the projectors 
\[ \pi_{0}(Y), \dots, \pi_{2 \dim(Y)}(Y) \]
provide a CK decomposition for $Y$.

Now we will prove the statement about the action. Recall the action of $\pi_{k}(Y)$ on $CH^{j}(Y)$ is given by the values of $\pi_{k}(Y)_{e}$, where $e= \Spec (k)$. We have to show that, for a given value of $p$, 
\[ \pi_{0}(Y)_{e}(\beta)= \cdots = \pi_{p-1}(Y)_{e}(\beta) = \pi_{2p+1}(Y)_{e}(\beta) = \dots = \pi_{2 \dim(Y)}(Y)_{e}(\beta)=0 \]
for any $\beta \in CH^{p}(Y)$.

As before, $\beta$ can be written as
\[ \sum_{r=0}^{p} \sum_{s=1}^{m_{r}} \pi^{\ast}(\alpha_{r,s}) \cap T_{r,s} \]
for some $\alpha_{r,s} \in CH^{p-r}(X)$.

From the calculations done before we see that
\[ \pi_{k}(Y)_{e}(\beta)= \sum_{\substack{(i,j) \in I_{k} \\ j \ even}} \dis \sum_{\lambda=1}^{m_{j/2}} \pi^{\ast}(\pi_{i}(X)_{e}(\alpha_{j/2, \lambda})) \cap T_{j/2, \lambda} \ , \]
and we have to consider two cases.

Suppose $0 \leq k \leq p-1$. Then for $(i,j) \in I_{k}$ we have
\[  0 \leq 2i+j \leq 2i +2j \leq 2(p-1) \]
and therefore
\begin{equation} \label{ineq01}
0 \leq i \leq p-\frac{j}{2}-1 \ .
\end{equation}

But $\alpha_{j/2, \lambda} \in CH^{p-\frac{j}{2}}(X)$ and $\pi_{i}(X)$ acts as zero there because of (\ref{ineq01}) and the hypothesis on $\pi_{i}(X)$. In this way we have that $\pi_{k}(Y)_{e}(\beta)=0$ for $0 \leq k \leq p-1$.

Now, suppose $2p+1 \leq k \leq 2\dim(Y)$. For $(i,j) \in I_{k}$, 
\[ 2p+1 \leq k= i+j \ ; \quad i \leq 2 \dim(X) \ . \]

Putting together these two inequalities we obtain $2(p-j/2)+1 \leq i \leq 2 \dim(X)$ and again, by the hypothesis on $\pi_{i}(X)$, it acts as zero on $CH^{p-j/2}(X)$. Therefore
$\pi_{k}(Y)_{e}(\beta)=0$ for $2p+1 \leq k \leq 2\dim(Y)$. \fin

\bibliographystyle{alpha}

\end{document}